\documentclass[reqno]{amsart}
\usepackage{amssymb}
\usepackage{amsmath}

\makeatletter
\@addtoreset{equation}{section}
\makeatother

\renewcommand\thefigure{\thesection.\@arabic\c@figure}
\renewcommand\thetable{\thesection.\@arabic\c@table}
 
\newtheorem{theorem}{Theorem}[section]
\newtheorem{lemma}[theorem]{Lemma}
\newtheorem{proposition}[theorem]{Proposition}

\newcommand{\mc}[1]{{\mathcal #1}}   
\newcommand{\mf}[1]{{\mathfrak #1}}
\newcommand{\mb}[1]{{\mathbf #1}}
\newcommand{\bb}[1]{{\mathbb #1}}

\newcommand{\<}{\langle}
\renewcommand{\>}{\rangle}
\newcommand{\eps}{\varepsilon}

\newcommand{\x}{{\bf x}}

\begin{document}
 
\author{M. D. Jara, C. Landim}

\address{\noindent IMPA, Estrada Dona Castorina 110,
CEP 22460 Rio de Janeiro, Brasil
\newline
e-mail:  \rm \texttt{monets@impa.br}
}

\address{\noindent IMPA, Estrada Dona Castorina 110,
CEP 22460, Rio de Janeiro, Brasil and CNRS UMR 6085,
Avenue de l'Universit\'e, BP.12, Technop\^ole du Madrillet,
F76801 Saint-\'Etienne-du-Rouvray, France. 
\newline
e-mail:  \rm \texttt{landim@impa.br}
}

\title[Quenched monequilibrium central limit theorem]{Quenched
  nonequilibrium central limit theorem for a tagged particle in the
  exclusion process with bond disorder}

\begin{abstract}  
  For a sequence of i.i.d. random variables $\{\xi_x : x\in \bb Z\}$
  bounded above and below by strictly positive finite constants,
  consider the nearest-neighbor one-dimensional simple exclusion
  process in which a particle at $x$ (resp. $x+1$) jumps to $x+1$
  (resp. $x$) at rate $\xi_x$. We examine a quenched nonequilibrium
  central limit theorem for the position of a tagged particle in the
  exclusion process with bond disorder $\{\xi_x : x\in \bb Z\}$. We
  prove that the position of the tagged particle converges under
  diffusive scaling to a Gaussian process if the other particles are
  initially distributed according to a Bernoulli product measure
  associated to a smooth profile $\rho_0:\bb R\to [0,1]$.
\end{abstract}



\thanks{C. Landim was partially supported by the John S. Guggenheim
  Memorial Foundation, FAPERJ and CNPq.}
 
\maketitle

\section{Introduction}
\label{sec0}

A classical problem in statistical mechanics consists in proving that
the dynamics of a single particle in a mechanical system is well
approximated on a large scale by a Brownian motion (\cite{s, kl}).  In
a seminal paper, Kipnis and Varadhan \cite{kv} proved an invariance
principle for the position of a tracer particle in the symmetric
simple exclusion process. The method relies on a central limit theorem
for additive functionals of Markov processes and uses time
reversibility and translation invariance. This approach has been
extended to interacting particle systems whose generators satisfy a
sector condition or, more generally, graded sector conditions
(\cite{lov} and references therein).

In \cite{JL}, we proved a non-equilibrium central limit theorem for
the position of a tagged particle in the one-dimensional
nearest-neighbor symmetric exclusion process. We assumed that the
initial state is a product measure associated to a smooth
profile. Observing that the position of the tagged particle can be
recovered from the density field and the total current through a bond,
we deduced a central limit theorem for the tagged particle from a
joint non-equilibrium central limit theorem for the density field and
the current.

The evolution of random walks in random environment has attracted some
attention in these last years (\cite{sz} and references therein).
Recently, a quenched central limit theorem has been proved for random
conductance models \cite{ss}. Here, to each bond $\{x,y\}$ of $\bb
Z^d$ is attached i.i.d.  strictly positive random variables
$\xi_{x,y}$. Under some conditions on the variables $\xi$, the authors
proved, among other results, that for almost all environments $\xi$, a
random walk on $\bb Z^d$ which jumps from $x$ to $y$ at rate
$\xi_{x,y}$ converges, when diffusively rescaled, to a Brownian
motion.

In this article we consider a one-dimensional nearest-neighbor
exclusion process evolving on an environment $\xi$. Each particle
behaves as the random walk described above, with the additional rule
that a jump is suppressed whenever a particle decides to jump over
a site already occupied. Under very mild assumptions on the
environment, we prove that the density field converges to the solution
of a heat equation, generalizing a previous result obtained by Nagy
\cite{Nag}. 

Assuming that the environment is strictly elliptic, i.e., formed by
i.i.d. random variables $\xi_{x,x+1}$ strictly bounded away from $0$
and $\infty$, we prove a non-equilibrium central limit theorem for the
density field, which holds for almost all realizations of the
environment. Here the assumption of independence and identical
distribution of the environment could be relaxed. In contrast with
\cite{gkl}, where annealed central limit theorems are considered, we
prove in this article a quenched statement.

From this result and from a non-equilibrium central limit theorem for
the current, we prove the main result of the article which states a
central limit theorem for the position of a tagged particle starting
from a configuration in which particles are distributed according to a
Bernoulli product measure associated to a smooth density profile. This
central limit theorem holds for almost all environment $\xi$'s.

The approach and the main technical difficulties can be summarized in
few words to the specialists. The model is in principle non-gradient
due to the presence of the environment \cite{fm}. However, a
functional transformation of the the empirical measure \eqref{f50}
turns it into a gradient model. The proof that the transformed
empirical measure is close to the original empirical measure imposes
some conditions on the the environment.

The same strategy can be applied to derive a nonequilibrium central
limit theorem for the density field. Here, however, to prove tightness
and to show that the transformed density field is close to the
original, some sharp estimates on the space time correlations are
needed as well. The deduction of these estimates require a Nash type
bound on the kernel of the random walk in the random conductance
model, which has been proved only under a strict ellipticity condition
of the environment. A this point, it remains to adapt the strategy
introduced in \cite{JL} to prove the central limit theorem for the
tagged particle.

While in Rome in April 2005, the second author showed to A.
Faggionato the model and the method described in next section to
derive the hydrodynamic behavior of this bond disorder model. At that
time he thought that the approach required uniform ellipticity of the
environment. A few months later and independently, Alessandra
\cite{af}, generalizing Nagy's method \cite{Nag}, obtained a proof of
the hydrodynamic behavior requiring only the assumptions stated in
Theorem \ref{t1} below, while the authors realized that strict
ellipticity was not needed in their approach.

\section{Main results}
\label{sec1}

We state in this section the main results of the article.  Denote by
$\mathcal X$ the state space $\{0,1\}^\mathbb Z$ and by the Greek
letter $\eta$ the elements of $\mc X$ so that $\eta(x)=1$ if there is
a particle at site $x$ for the configuration $\eta$ and $\eta(x)=0$
otherwise.

Consider a sequence $\{\xi_x: x\in \bb Z\}$ of strictly positive
numbers.  The symmetric nearest-neighbor simple exclusion process with
bond disorder $\{\xi_x: x\in \bb Z\}$ is the Markov process $\{\eta_t
: t \geq 0\}$ on $\mathcal X$ whose generator $L_{\xi, N}$ acts on
cylinder functions $f$ as $$
(L_{\xi, N} f) (\eta) = N^2 \sum_{x \in
  \mathbb Z} \xi_x \, (\nabla_x f) (\eta), $$
where $(\nabla_x f)
(\eta)=f(\sigma^{x,x+1}\eta)-f(\eta)$ and $$
\sigma^{x,y}\eta(z)=
\begin{cases}
\eta(y), & z=x\\
\eta(x), & z=y\\
\eta(z), & z\neq x,y.\\
\end{cases}
$$
Notice that the process is speeded up by $N^2$.

Existence and ergodic properties of this Markov process can be proved
as in the space homogeneous case (\cite{Lig},\cite{Nag}). Moreover,
the Bernoulli product measures $\nu_\alpha$ in $\{0,1\}^\mathbb Z$,
with marginals $\nu_\alpha\{\eta(x)=1\}=\alpha$ for $\alpha \in
[0,1]$, are extremal, reversible measures.  

For each profile $\rho_0: \mathbb R \to [0,1]$, denote by
$\nu^N_{\rho_0(\cdot)}$ the product measure on $\mc X$ with marginals
given by $\nu^N_{\rho_0(\cdot)}\{ \eta (x)=1\} =\rho_0(x/N)$.  For a
measure $\mu$ on $\mc X$, let $\bb P_{\mu}^N$ stand for the
probability measure on the path space $D(\bb R_+, \mc X)$ induced by
the Markov process $\eta_t$ and the measure $\mu$.

The empirical measure associated to the process $\eta_t$ is defined
by
$$
\pi^N_t(du) = \frac{1}{N}\sum_{x \in \mathbb Z} 
\eta_{t}(x) \delta_{x/N}(du).
$$

Fix $0<\gamma<\infty$.  Let $C_{0}^{2}(\mathbb R)$ be the set of twice
continuously differentiable functions $G: \mathbb R \to \mathbb R$
with compact support. Fix a profile $\rho_0 : \bb R_+ \to [0,1]$.  A
bounded function $\rho : \bb R_+ \times \bb R \to [0,1]$ is said to be
a weak solution of the heat equation
\begin{equation}
\label{f8}
\left\{
\begin{array}{l}
\partial_t \rho = \gamma^{-1} \Delta \rho \\
\rho(0,\cdot) =\rho_0(\cdot)\end{array}
\right.
\end{equation}
if 
\begin{equation}
\label{f7}
\< \rho_t, G \> \;=\; \< \rho_0, G \> \;+\;
\int_0^t ds\, \< \rho_s, \gamma^{-1} \Delta G \> 
\end{equation} 
for all $t\ge 0$ and all $G$ in $C_{0}^{2}(\mathbb R)$. In these
equations, $\Delta$ stands for the Laplacian and $\< \rho , H \>$ for
the integral of $H$ with respect to the measure $\rho(u) du$. It is
well known that for any bounded profile $\rho_0: \bb R\to [0,1]$,
there exists a unique weak solution of \eqref{f8}. The first main
result of the article states a quenched law of large numbers for the
empirical measure under weak assumptions on the environment $\{\xi_x :
x\in \bb Z\}$.

\begin{theorem}
\label{t1}
Assume that
\begin{equation}
\label{f9}
\lim_{K\to\infty} \frac 1K \sum_{x=1}^K \xi_x^{-1} \;=\; \gamma\; ,
\quad
\lim_{K\to\infty} \frac 1K \sum_{x=-K}^{-1} \xi_x^{-1} \;=\; \gamma
\end{equation}
for some $0< \gamma < \infty$. Fix a profile $\rho_0 : \bb R \to
[0,1]$. Under $\bb P_{\nu^N_{\rho_0(\cdot)}}^N$, $\pi^N_{t}$ converges
in probability to the weak solution of \eqref{f8}: For every
continuous function with compact support $G$, every $t\ge 0$ and every
$\delta >0$,
\begin{equation*}
\lim_{N\to\infty} \bb P_{\nu^N_{\rho_0(\cdot)}}^N \Big[ \, \big\vert
\< \pi^N_t, G\> - \< \rho_t, G\> \big\vert > \delta \Big] \;=\; 0\;,
\end{equation*}
where $\rho(t,u)$ is the weak solution of \eqref{f8}.
\end{theorem}

To prove a quenched nonequilibrium central limit theorem for the
empirical measure, assume that $\{\xi_x : x\in \bb Z\}$ is a sequence
of i.i.d.  random variables defined on a probability space $(\Omega,
P, \mc F)$ such that 
\begin{equation}
\label{f31}
P[\varepsilon \le \xi_0 \le \varepsilon^{-1}]=1
\end{equation}
for some $\varepsilon>0$. This strong ellipticity condition is needed
in Section \ref{sec4} to prove sharp estimates of the decay of the
space-time correlation functions. All other arguments require the
weaker integrability condition: $E[\xi^{-6}_0]<\infty$.

Fix a profile $\rho_0:\bb R\to [0,1]$ and an environment $\xi$.  Let
$\rho^{N,\xi}_t(x) = \bb E_{\nu^N_{\rho_0(\cdot)}} [\eta_t (x)]$. A
trivial computation shows that $\rho^{N,\xi}_t:\bb Z\to [0,1]$ is the
solution of the discrete linear equation
\begin{equation}
\label{f16}
\left\{
\begin{array}{l}
\partial_t \rho_t (x) \;=\; N \big\{ \xi_x (\nabla_N \rho_t) (x)
- \xi_{x-1} (\nabla_N \rho_t) (x-1) \big\} \;, \\
\rho_0 (x) \;=\; \rho_0(x/N)\;,
\end{array}
\right.
\end{equation}
where $(\nabla_N h)(x) = N\{h(x+1) - h(x)\}$. We denote frequently
$\rho^{N,\xi}_t$ by $\rho^{N}_t$.

Denote by $\mathcal S(\mathbb R)$ the Schwartz space of rapidly
decreasing functions and by $\mathcal S'(\mathbb R)$ its dual, the
space of distributions. Let $\{Y^N_t, t\ge 0\}$ be the density
fluctuation field, a $\mathcal S'(\mathbb R)$-valued process given by
$$
Y_t^N (G) \;=\; \frac{1}{\sqrt{N}} \sum_{x \in \bb Z}
G(x/N) \{ \eta_t (x) - \rho^{N,\xi}_t(x)\}
$$
for $G$ in $\mc S(\mathbb R)$. Next theorem states the almost sure
convergence of the finite dimensional distributions of $Y^N_t$ to the
marginals of a centered Gaussian field.

\begin{theorem}
\label{t2}
Let $\{\xi_x : x\in \bb Z\}$ be a sequence of i.i.d. random variables
satisfying assumption \eqref{f31}.  Let $\rho_0: \bb R \to [0,1]$ be a
profile with first derivative in $L^1(\bb R) \cap L^\infty(\bb R)$.
There exists a set of environments $\Omega_0$ with total measure such
that for every $\xi$ in $\Omega_0$, every $k\ge 1$ and every $0\le t_1
< \cdots < t_k$, $(Y^N_{t_1}, \dots, Y^N_{t_k})$ converges to a
centered Gaussian vector $(Y_{t_1}, \dots, Y_{t_k})$ with covariance
given by
\begin{eqnarray}
\label{f17}
E[Y_s(G) Y_t(H)] &\!\!\!=\!\!\!&  \int_{\bb R}
\chi(\rho_0(u)) \, T_{s} G (u) \, T_{t} H (u) \\
&\!\!\!+\!\!\!& 2 \gamma^{-1} \int_0^s dr\, \int_{\bb R}
\chi(\rho(r,u)) \, \nabla T_{s-r} G (u) \, \nabla T_{t-r} H (u) 
\nonumber
\end{eqnarray}
for all $0\le s\le t$, $G$, $H$ in $\mc S(\bb R)$. Here $\rho$ stands
for the solution of the heat equation \eqref{f8}, $\{T_r : r\ge 0\}$
for the semi-group associated to $\gamma^{-1} \Delta$ and
$\chi(\alpha) = \alpha (1-\alpha)$ for the compressibility in the
exclusion process.
\end{theorem}

Denote by $\nu^{N,*}_{\rho_0(\cdot)}$ the measure
$\nu^N_{\rho_0(\cdot)}$ conditioned to have a particle at the origin
and by $X_t^N$ the position at time $t$ of the particle initially
at the origin. Define $u_t^N = u_t^{N,\xi}$ by the relation
\begin{equation}
\label{f32}
\sum_{x=0}^{u_t^N} \rho_t^{N,*}(x) \; \leq \; \xi_{-1}
\int_0^t  N^2 \{ \rho^{N,*}_s(-1) - \rho^{N,*}_s(0) \} \, ds 
\; <\; \sum_{x=0}^{u_t^N+1} \rho_t^{N,*}(x)\;,
\end{equation}
where $\rho^{N,*}_t$ is the solution of \eqref{f16} with initial
condition $\rho^{N,*}_0(0) = 1$, $\rho^{N,*}_0(x) = \rho_0(x/N)$,
$x\not = 0$.  Let $W_t^N = (X_t^N - u_t^N)/\sqrt N$.

\begin{theorem}
\label{t3}
Let $\{\xi_x : x\in \bb Z\}$ be a sequence of i.i.d. random variables
satisfying assumption \eqref{f31}.  Let $\rho_0$ be an initial profile
with first derivative in $L^1(\bb R) \cap L^\infty(\bb R)$ and second
derivative in $L^\infty(\bb R)$. There exists a set of environments
$\Omega_0$ with total measure and the following property. For every
$\xi$ in $\Omega_0$, every $k \geq 1$ and every $0 \leq t_1<...<t_k$,
under $\bb P_{\nu_{\rho_0(\cdot)}^{N,\ast}}$,
$(W^N_{t_1},...,W^N_{t_k})$ converges in law to a Gaussian vector
$(W_{t_1},...,W_{t_k})$ with covariances given by
\begin{eqnarray*}
\rho(s,u_s) \rho(t,u_t) E[W_s W_t] &=&
\int_{-\infty}^0 dv \,  P[Z_s\le v]
\, P[Z_t\le v] \, \chi(\rho_0(v)) \\
&+& \int_0^{\infty} dv \,  P [Z_s\ge v] \, P[Z_t\ge v]
\, \chi(\rho_0(v))  \\
&+& \frac 2\gamma \int_0^s dr \int_{-\infty}^{\infty} dv 
\, p_{t-r}(u_t,v) \, p_{s-r}(u_s,v) \, \chi(\rho(r,v))
\end{eqnarray*}
provided $s\le t$. In this formula, $Z_t = u_t + B^0_{t/\gamma}$, where
$B^0_t$ is a standard Brownian motion starting from the origin, and
$p_t(v,w)$ stands for the kernel of $B^0_{t/\gamma}$.
\end{theorem}

\section{Hydrodynamic limit}
\label{sec2}

We prove in this section Theorem \ref{t1}. Fix an environment
satisfying \eqref{f9} and denote by $\mc M_+(\mathbb R)$ the set of
positive Radon measures in $\mathbb R$. Fix $T\ge 0$ and a bounded
profile $\rho_0 : \bb R\to [0,1]$.  Let $\{Q_N : N \ge 1\} =
\{Q_{N,\xi} : N \ge 1\}$ be the sequence of measures on $D([0,T],
\mathcal M_+(\mathbb R))$ induced by the Markov process $\pi^N_t$ and
the initial state $\nu^N_{\rho_0(\cdot)}$.

The proof of Theorem \ref{t1} is divided in two steps.  We first prove
tightness of the sequence $\{Q_N\}_N$, and then that all limit points
of $\{Q_N\}_N$ are supported on weak solutions of the hydrodynamic
equation. It follows from these two results and the uniqueness of weak
solutions of the heat equation \eqref{f8} that $\pi^N_t$ converges in
probability to the absolutely continuous measure $\rho(t,u) du$ whose
density is the solution of \eqref{f8} (cf. \cite{kl}).

It turns out that this program can not be accomplished for the
empirical measure $\pi^N_t$, but for a ``corrected by the
environment'' process $X^N_t$, which is close enough to the empirical
measure $\pi^N_t$.

\subsection{Corrected empirical measure}

Denote by $C^2_0(\bb R)$ the space of twice continuously
differentiable functions with compact support. For a function $G$ in
$C^2_0(\bb R)$ and an environment $\xi$, let $T_\xi G : \bb Z\to\bb R$
be the sequence defined by
\begin{equation}
\label{f1}
(T_\xi G)(x) \;=\; \sum_{j < x} \xi^{-1}_j \Big\{
G\big(\frac{j+1}{N}\big)-G\big(\frac{j}{N}\big) \Big\}\;.
\end{equation}

For each each $N\ge 1$ and each function $G$ in $C^2_0(\bb R)$, the
series $\sum_x \xi_x^{-1}[G((x+1)/N) -G(x/N)]$ is absolutely summable
because $G$ has compact support. Moreover, it follows from \eqref{f9}
that 
\begin{equation}
\label{f6}
T_{\xi, G} \;=\; T^N_{\xi, G} \;:=\; \sum_{x\in\bb Z} \xi_x^{-1} 
\big\{ G((x+1)/N) -G(x/N) \big\}
\end{equation}
converges to $0$ as $N \uparrow \infty$.

We introduce $T_\xi G$ for two reasons. On the one hand, we expect
$(T_\xi G)(x)$ to be close to $ \gamma G(x/N)$, which is the content
of Lemma \ref{s2} below. On the other hand,
\begin{equation*}
N \big\{ (T_\xi G)(x+1) - (T_\xi G)(x) \big\} \, \xi_x \;=\; (\nabla_N
G)(x/N)\; ,
\end{equation*}
where $\nabla_N$ stands for the discrete derivative: $(\nabla_N
G)(x/N) = N\{ G(x+1/N) - G(x/N) \}$. Hence, formally,
\begin{equation}
\label{f4}
L_{\xi,N} \frac 1N \sum_{x\in\bb Z} (T_\xi G)(x) \eta(x) \;=\;
\frac 1N \sum_{x\in\bb Z} (\Delta_N G)(x/N) \eta(x)\;,
\end{equation}
where $\Delta_N$ stands for the discrete Laplacian. 

Of course, $T_\xi G$ may not belong to $\ell_1(\bb Z)$, the space of
summable series, and the left hand side of the previous formula may
not be defined. To overcome this difficulty, we modify $T_\xi G$ in
order to integrate it with respect to the empirical measure.  Fix an
arbitrary integer $l > 0$ which will remained fixed in this section.
let $g=g_l:\mathbb R \to \mathbb R$ be defined by
$$
g (u)=
\begin{cases}
0, &u < 0\;, \\
u/l, &0 \leq u <l\;, \\
1, &u \geq l\;.
\end{cases}
$$
For each function $G$ in $C^2_0(\bb R)$, let
\begin{equation}
\label{f2}
(T_{\xi,l} G)(x) \;:=\; (T_{\xi} G)(x) \;-\; \frac 
{T_{\xi, G}} {T_{\xi, g}} (T_{\xi} g)(x)\; .
\end{equation}
Notice that $T_{\xi, g}$ converges to $\gamma$ almost surely, as
$N\uparrow\infty$. In particular, by \eqref{f6} the ratio $T_{\xi,
  G}/T_{\xi, g}$ vanishes almost surely as $N\uparrow\infty$.  In the
end of this section we prove the following statement.

\begin{lemma}
\label{s2}
For each function $G$ in $C^2_0(\bb R)$, and each environment $\xi$
satisfying \eqref{f9}, $T_{\xi,l} G$ belongs to $\ell_1 (\bb Z)$ and
\begin{equation*}
\lim_{N\to\infty} \frac{1}{N} \sum_{x \in \mathbb Z} 
\big| T_{\xi,l} G (x) - \gamma \, G(x/N) \big| \;=\; 0\; .
\end{equation*}
\end{lemma}

Denote by $X^N_t$ the corrected empirical measure defined by
\begin{equation}
\label{f50}
X^N_t(G) \;=\; X^{N, l, \xi}_t(G) \;=\; 
\frac {1}{N}\sum_{x \in \mathbb Z} T_{\xi,l} G (x) \, 
\eta_t^N(x)
\end{equation}
for each function $G$ in $C^2_0(\bb R)$. 

As mentioned before, the sequence $T_{\xi,l} G (x)$ has two
properties. On the one hand, in view of Lemma \ref{s2}, it is close to
$\gamma G(x/N)$ in $\ell_1(\bb Z)$. In particular, the integral of $G$
with respect to the empirical measure is close to $\gamma^{-1}
X^N_t(G)$ uniformly in time.  On the other, by \eqref{f4}, the
martingale associated to $\gamma^{-1} X^N_t(G)$ has an integral term
which can be expressed as function of the empirical measure. Indeed,
for a function $G$ in $C^2_0(\bb R)$, let $M^N_t(G) = M^{N, l,
  \xi}_t(G)$ be the martingale defined by
\begin{eqnarray}
\label{f3}
\!\!\!\!\!\!\!\!\!\!\!\! &&
M^N_t(G) \;=\; X^N_t(G) \;-\; X^N_0(G) \;-\; \int_0^t ds\,
N^2 L X^N_s(G) \\
\!\!\!\!\!\!\!\!\!\!\!\! && \quad
\;=\; X^N_t(G) \;-\; X^N_0(G) \;-\; \int_0^t ds\,
\Big\{ \< \pi^N_s , \Delta_N G \> - \frac 
{T_{\xi, G}} {T_{\xi, g}} \< \pi^N_s , \Delta_N g \> \Big\} \;.
\nonumber
\end{eqnarray}

\subsection{Tightness of $\pi_t^N$.}

It is well known that a sequence of probability measures $\{Q_N\}_N$
on $D([0,T], \mc M_+(\bb R))$ is tight if and only if the sequence
$\{Q_N(G)\}_N$ is tight for all $G \in C^2_0(\bb R)$, where $Q_N(G)$
is the probability measure in $D([0,T],\mathbb R)$ corresponding to
the process $\< \pi^N_t, G\>$.

We claim that the process $X^N_t(G)$ is tight.  Recall Aldous
criteria for tightness in $D([0,T],\mathbb R)$:

\begin{lemma}
\label{s4}
A sequence of probability measures $\{P_N\}_N$ in $\mathcal
D([0,T],\mathbb R)$ is tight if
\begin{itemize}
\item[(i)] For all $0\le t\le T$ and for all $\eps > 0$ there exists a
  finite constant $A$ such that  $\sup_N P_N(|x_t|>A) < \eps$, 
\item [(ii)] For all $\delta > 0$, 
\begin{equation*}
\lim_{\beta \to 0} \limsup_{N \to \infty}
\sup_{\substack{\tau \in \mathcal T\\\theta \leq \beta}}
P_N(|x_{\tau + \theta}-x_\tau|>\delta)=0\;,
\end{equation*}
where $\mathcal T$ is the set of stopping times with respect to the
canonical filtration bounded by $T$.
\end{itemize}
\end{lemma}

To prove tightness of $X^N_t(G)$ note that (i) is automatically
satisfied because the number of particles per site is bounded and
$T_{\xi,l} G$ converges to $\gamma G$ in $\mathcal \ell_1(\mathbb Z)$.

To check condition (ii), fix a stopping time $\tau$ bounded by $T$ and
$\theta \leq \beta$. Recall from formula \eqref{f3} that we may
express $X^N_{\tau + \theta} (G) - X^N_{\tau} (G)$ as the sum of a
martingale difference and an integral.  On the one hand, computing the
quadratic variation of the martingale $M_{t}^N(G)$, we obtain that
\begin{eqnarray*}
\!\!\!\!\!\!\!\!\!\!\!\!\!\! &&
\mathbb E_N[(M_{\tau + \theta}^N(G) - M_{\tau}^N(G))^2] \; = \\
\!\!\!\!\!\!\!\!\!\!\!\!\!\! && \quad
\mathbb E_N \Big[\int_\tau^{\tau +\theta} ds\, 
\frac{1}{N^2} \sum_{x \in \mathbb Z} 
\Big \{\nabla_N G(x/N) - \frac{T_{\xi, G}}{T_{\xi, g}} \nabla_N g(x/N)
\Big\}^2 (\eta_s^N(x+1)-\eta_s^N(x))^2 \Big]\; .
\end{eqnarray*}
The previous expression is bounded above
by $N^{-1}\{ C(G) + l^{-1} (T_{\xi, G}/T_{\xi, g})^2\}$, which
vanishes as $N\uparrow\infty$ in view of \eqref{f9}.

On the other hand, since there is at most one particle per site and
since $G$ belongs to $C^2_0(\bb R)$,
\begin{equation*}
\Big|\int_\tau^{\tau + \theta} ds \, \frac{1}{N} \sum_{x \in \mathbb
Z} \Big \{\Delta_N G(x/N)-  \frac{T_{\xi, G}}{T_{\xi, g}} 
\Delta_N g (x/N) \Big\} \eta_s^N(x) \, \Big| \; \leq \;
C_0 \beta  \;+\;  \frac{2 \beta}{l} 
\Big\vert \frac{T_{\xi, G}}{T_{\xi, g}}\Big\vert
\end{equation*}
for some finite constant $C_0$ depending only on $\beta$. As
$N\uparrow\infty$, the second term vanishes in view of \eqref{f9},
\eqref{f6}.  This proves condition (ii) of Lemma \ref{s4} and
tightness of the process $X^N_t(G)$.

In view of Lemma \ref{s2},
\begin{equation}
\label{f5}
\sup_{0\le t\le T} |\gamma \< \pi_t^N,G \> - X^N_t(G) | 
\;\leq\; \frac{1}{N} \sum_{x \in \mathbb Z}
|\gamma G(x/N)- T_{\xi,l} G (x)|
\end{equation}
converges to $0$ as $N \uparrow 0$. In particular, $\< \pi_t^N,G \>$ is
also tight, with the same limit points of $X^N_t(G)$. Since this
statement holds for all $G$ in $C^2_0(\bb R)$, the sequence $Q_N$ is
tight.

\subsection{Uniqueness of limit points}

Let $Q$ be a limit point of the sequence $\{Q_N\}_N$. Since there is
at most one particle per site, $Q$ is concentrated on absolutely
continuous paths $\pi(t,du) = \rho(t,u) du$, with positive density
bounded by $1$: $0\le \rho(t,u) \le 1$.

We have seen in the last subsection that $Q$ is also a limit point of
$X^N_t$. Fix a function $G$ in $C^2_0(\bb R)$ and recall the
definition of the martingale $M^N_t(G)$ given in \eqref{f3}.  By the
proof of the tightness of $X^N_t$, the expectation of the quadratic
variation of $M^N_t(G)$ vanishes as $N\uparrow\infty$.  In particular,
in view of \eqref{f5}, the measure $Q$ is concentrated in trajectories
$\pi_t$ such that 
\begin{equation*}
\< \pi_t, G \> \;=\; \< \pi_0, G \> \;+\;
\int_0^t ds\, \< \pi_s, \gamma^{-1} \Delta G \> 
\end{equation*}
for all $0\le t\le T$, $G$ in $C_0^2(\bb R)$. By the uniqueness of
weak solutions of the heat equation, Theorem \ref{t1} is proved. \qed
\medskip

We conclude this section with the 

\noindent{\bf Proof of Lemma \ref{s2}.} $T_{\xi,l} G$
belongs to $\ell_1(\bb Z)$ because it belongs to $\ell_\infty(\bb Z)$
and vanishes outside a finite set.  Fix a smooth function $G$ in
$C^2_0(\bb R)$.

Consider first the sum over $x\le 0$. In this case $(T_\xi g_l)(x) =0$
so that $(T_{\xi,l} G)(x) = (T_{\xi} G)(x)$. In particular,
\begin{equation*}
\frac 1N \sum_{x\le 0} \Big\vert (T_{\xi,l} G)(x) - \gamma\, G(x/N) 
\Big\vert \;=\; \frac 1N \sum_{x\le 0} \Big\vert
\frac 1N \sum_{y< x} \hat \xi_y^{-1} (\nabla_N G)(y/N) \Big\vert\; , 
\end{equation*}
where $\hat \xi_y^{-1} = \xi_y^{-1} - \gamma$. Both sums in $x$ and
$y$ start from $-AN$, for some $A>0$, because $G$ has compact support.
Fix $\varepsilon >0$. Since $G'$ is uniformly continuous, there exists
$\delta >0$ such that $|G'(v) - G'(u)|\le \varepsilon$ if $|v-u|\le
\delta$. We may therefore replace $(\nabla_N G)(y/N)$ by
$G'(k\delta)$, for $k\delta \le y/N \le (k+1)\delta$, paying a price
bounded by $C(G) \varepsilon$. After this replacement, the law of
large numbers \eqref{f9} ensures that the previous expression
vanishes as $N\uparrow\infty$.

Similarly, for $x\ge l N$, $(T_\xi g_l)(x) = T_{\xi,g_l}$ so that
$(T_{\xi,l} G)(x) = (T_{\xi} G)(x) - T_{\xi,G}$. Therefore,
\begin{eqnarray*}
(T_{\xi,l} G)(x) - \gamma\, G(x/N) \;=\; 
- \frac 1N \sum_{y\ge x} \hat \xi_y^{-1}  (\nabla_N G)(y/N) 
\end{eqnarray*}
and we may repeat the previous arguments to show that the sum for
$x\ge lN$ vanishes as $N\uparrow\infty$.
 
Finally, for $0\le x < lN$, we estimate separately $(T_{\xi} G)(x) -
\gamma G(x/N)$ and $\{T_{\xi,G}/$ $T_{\xi,g_l}\} (T_{\xi} g_l)(x)$.
The first piece is handled as before, while the second vanishes as
$N\uparrow\infty$ in view of \eqref{f6} and because $(T_{\xi} g_l)(x)
/ T_{\xi,g_l}$ is absolutely bounded by $1$. This proves Lemma
\ref{s2}. \qed

\section{Fluctuations of the empirical measure}
\label{sec3}

Let $\{\xi_x : x\in \bb Z\}$ be a sequence of i.i.d. random variables
defined on a probability space $(\Omega, P, \mc F)$. We prove in this
section a quenched nonequilibrium central limit theorem for the
empirical measure. The proof relies on sharp estimates of the decay of
the space-time correlation functions presented in Section \ref{sec4}
which requires the strong ellipticity condition: $P[\varepsilon \le
\xi^{-1}_0 \le \varepsilon^{-1}]=1$ for some $\varepsilon>0$. To
stress that it is only in the estimation of the correlation functions
that we need this condition, we present all other proofs under the
weaker assumption that $E[\xi_0^{-6}]<\infty$.  Moreover, the
hypotheses of independence, identical distribution and finiteness of
the sixth moment can be relaxed. 

Throughout this section the index $l$ of the operator $T_{\xi,l}$
introduced in the previous section depends on $N$ as $l=l_N=N^{1/4}$.
Recall that we denote by $\mc S(\bb R)$ the Schwartz space of rapidly
decreasing functions. We may extend the operators $T_\xi$, $T_{\xi,l}$
to $\mc S(\bb R)$:

\begin{lemma}
\label{s3}
Assume that $E[\xi_0^{-6}]<\infty$ and fix a function $G \in \mc S(\bb
R)$. There exists a subset $\Omega_G$ with total measure such that for
each $\xi$ in $\Omega_G$ $T_\xi G (x)$ is well defined and
\begin{equation*}
\lim_{N\to\infty} N^{1/4} \sup_{x \in \bb Z} 
|T_\xi G(x) -\gamma G(x/N)| \;=\; 0\;.
\end{equation*}
In particular, $\lim_{N\to\infty} N^{1/4} T_{\xi,G} =0$. Moreover,
\begin{eqnarray*}
\!\!\!\!\!\!\!\!\!\!\!\!\!\! &&
\lim_{N\to\infty} \sup_{x \in \bb Z} 
|T_{\xi,l} G(x) -\gamma G(x/N)| \;=\; 0 \\
\!\!\!\!\!\!\!\!\!\!\!\!\!\! && \quad \text{and}\quad
\lim_{N\to\infty} \frac{1}{N} \sum_{x \in \bb Z} 
|T_{\xi,l} G(x) -\gamma G(x/N)| \;=\; 0\;.
\end{eqnarray*}
\end{lemma}

The proof of this lemma is given at then end of this section.  By
interpolation it follows from this result that
\begin{equation}
\label{f15}
\lim_{N\to\infty} \frac{1}{N} \sum_{x \in \bb Z} |T_{\xi,l} G(x)|^{p} 
\;=\; \int |\gamma G(u)|^{p} \, du
\end{equation}
$\xi$-almost surely for all $1\le p\le \infty$.

Recall that the definition of the density field $Y^N_t$ given just
before the statement of Theorem \ref{t2}.  Denote by $Z_t^N$ the
fluctuation density field corrected by the environment: $$
Z_t^N(G)
\;=:\; \frac 1\gamma Y_t^N(T_{\xi,l}G) \;=\; \frac{1}
{\gamma \sqrt N}\sum_{x \in \mathbb Z}
(T_{\xi,l}G)(x) \{ \eta_t^N (x) - \rho^{N,\xi}_t(x)\} $$
for functions $G$ in $\mc S(\bb R)$.

We prove in this section a nonequilibrium central limit theorem for
the density field $Z^N_t$ in random environment and deduce from this
result the convergence of the finite dimensional distributions of the
field $Y^N_t$ defined in Section \ref{sec1}. Recall that we denote by
$\mc S'(\bb R)$ the Schwartz space of distributions. For a profile
$\rho_0 : \bb R \to (0,1)$ and an environment $\xi = \{\xi_x : x\in\bb
Z\}$ and $T>0$, let $Q^{N,\xi}_{\rho_0}$ be the measure on $D([0,T],
\mc S'(\bb R))$ induced by the process $Z^N_t$ and the initial state
$\nu_{\rho_0(\cdot)}^N$.

\begin{proposition}
\label{s1}
Fix a profile $\rho_0 : \bb R \to (0,1)$ with bounded and integrable
first derivative. There exists a set of environments $\Omega_0$ with
total measure such that for each $\xi$ in $\Omega_0$,
$Q^{N,\xi}_{\rho_0}$ converges to a centered Gaussian field $Z_t$ with
covariance given by \eqref{f17}.
\end{proposition}

The strategy of the proof of Proposition \ref{s1} is similar to the
one adopted for the hydrodynamic limit.  We prove tightness of the
distributions of $Z_t^N$ in $ D([0,T],\mathcal S'(\mathbb R))$ and
that all limit points of $Z_t^N$ satisfies a martingale problem which
characterizes the limiting measure.

We start proving tightness. For a function $G$ in $\mc S(\bb R)$,
consider the martingale $M_t^N(G)$ defined by
\begin{equation}
\label{f14}
M_t^N(G) \;=\; Z_t^N(G) \;-\; Z_0^N(G) \;-\; \int_0^t 
\gamma_1^N (s,G) \, ds\;,
\end{equation}
where
\begin{equation*}
\gamma_1^N (s,G)\;=\; Y_s^N( \gamma^{-1} \Delta_N G) 
\;-\; \frac{T_{\xi,G}} {T_{\xi,g}} \,  Y_s^N ( \gamma^{-1} \Delta_N g_l) \;.
\end{equation*}
The quadratic variation $\< M^N(G)\>_t$ of this martingale is equal to
$\int_0^t \gamma_2^N (s,G) \, ds$, with $\gamma_2^N (s,G)$ given
by
\begin{equation*}
\frac{1}{\gamma^2 N} \sum_{x \in \mathbb Z} \xi_x^{-1} 
\Big\{ (\nabla_N G) (x/N) - \frac{T_{\xi,G}}{T_{\xi,g}} 
(\nabla_N g_l) (x/N) \Big\}^2  (\eta_s^N(x+1)-\eta_s^N(x))^2\;.
\end{equation*}
In view of Mitoma's criterion for the relative compactness of a
sequence of measures in $D([0,T], \mathcal S'(\mathbb R))$ (cf.
\cite{Mit}, \cite{HS}, \cite{FPSV}), to show that the process $Z^N_t$
is tight, it is enough to prove that
\begin{equation}
\label{f10}
\sup_N \sup_{0 \leq t \leq T} \bb E_{\nu_{\rho_0(\cdot)}^N} 
[Z_t^N (G)^2] \;<\; \infty\;, \quad
\sup_N \sup_{0 \leq t \leq T} \bb E_{\nu_{\rho_0(\cdot)}^N}
[\gamma_i^N (t,G)^2] \;<\; \infty
\end{equation}
for $i=1$, $2$ and a dense family of functions $G$ in $\mc S(\mathbb
R)$.  Moreover, to show that all limit points of the sequence $Z^N_t$
are concentrated on $C([0,T],\mathcal S'(\mathbb R))$, we need to
check that for each function $G$ in $\mc S(\mathbb R)$ there exists 
a sequence $\delta_N =\delta(t,G,N)$, vanishing as $N\uparrow\infty$,
such that
\begin{equation}
\label{f11}
\lim_{N \to \infty} \bb P_{\nu_{\rho_0(\cdot)}^N} \Big [\sup_{0 \leq t \leq T}
|Z_t^N(G)-Z_{t-}^N(G)| \geq \delta_N \Big] \;=\; 0  \;.
\end{equation}
 
To prove \eqref{f10}, consider a countable dense subset of functions
$\mc S_0 (\bb R) = \{G_k : k\ge 1\}$ in $\mc S (\bb R)$. Let $\Omega_0
= \bigcap_{k\ge 1} \{ \Omega_{G_k} \cap \Omega_{( G_k')^2} \cap
\Omega_{\Delta G_k} \}$, where $\Omega_G$ are the total measure sets
introduced in Lemma \ref{s3}. Fix an environment $\xi$ in $\Omega_0$
and a function $G$ in the class $\{G_k : k\ge 1\}$. By Theorem
\ref{s11},
\begin{equation*}
\mathbb E_{\nu_{\rho_0(\cdot)}^N} [Z_t^N(G)^2] \;\le\;
\frac{1}{\gamma^2 N} \sum_{x \in \mathbb Z} (T_{\xi,l}G)(x)^2 \;+\;
\frac{C_1}{\gamma^2} \sqrt T \Big( \frac{1}{N} \sum_{x \in \mathbb Z} 
|(T_{\xi,l}G)(x)| \Big)^2
\end{equation*}
for some finite constant $C_1$ depending only on $\varepsilon$ and
$\rho_0$.  Since $\xi$ belongs to $\Omega_0$, by \eqref{f15}, as
$N\uparrow\infty$, these expressions converge to finite expressions.

On the other hand, by definition of $\gamma_1^N(t,G)$, 
$$
\frac{\gamma^2}{2} \mathbb E_{\nu_{\rho_0(\cdot)}^N}
[\gamma_1^N(t,G)^2] \;\leq\; \mathbb E_{\nu_{\rho_0(\cdot)}^N}
[Y_t^N(\Delta_N G)^2] \;+\; \Big\{\frac{T_{\xi,G} }{ T_{\xi,g}}
\Big\}^2 E_{\nu_{\rho_0(\cdot)}^N} [Y_t^N(\Delta_N g_l)^2] \;.  
$$
The first term is handled in the same way as $Z^N_t(G)$. The second 
term is also simple to estimate, since 
$$
Y_t^N(\Delta_N g_l)^2 \;=\;
\frac{N}{l_N^2} \Big\{ \big[ \eta_t^N (0) - \rho_t^N (0) \big]
- \big[ \eta_t^N(lN) - \rho_t^N(lN) \big] \Big\}^2 \;, 
$$
and since, by Lemma \ref{s3}, $(T_{\xi,G} / T_{\xi,g})^2 N
l_N^{-2}$ vanishes as $N \uparrow \infty$ for all $\xi$ in $\Omega_0$.

Finally, by definition of $\gamma_2^N(t,G)$,
\begin{equation*}
\mathbb E_{\nu_{\rho_0(\cdot)}^N} [\gamma_2^N(t,G)^2] 
\;\leq\; \frac{2}{\gamma^2 N} \sum_{x \in \mathbb Z} \xi_x^{-1} 
\nabla_N G(x/N)^2
\;+\; \frac{2}{N l_N^2 \gamma^2} \Big(\frac{T_{\xi,G}}{ T_{\xi,g}}\Big)^2 
\sum_{0\leq x < lN} \xi_x^{-1}\;.
\end{equation*}
The first term converges to a finite constant as $N \uparrow \infty$,
while the second term vanishes for $\xi$ in $\Omega_0$.

Since condition \eqref{f11} follows from the fact that no more than
one particle jumps at each time, the previous estimates show that for
each environment $\xi$ in $\Omega_0$, the sequence
$Q^{N,\xi}_{\rho_0}$ is tight and that each limiting point is
concentrated on $C([0,T],\mathcal S'(\mathbb R))$.  \medskip

We consider now the question of uniqueness of limit points.  Fix $\xi$
in $\Omega_0$, let $Q^\xi$ be a limit point of $Q^{N,\xi}_{\rho_0}$
and assume without loss of generality that $Q^{N,\xi}_{\rho_0}$
converges to $Q^\xi$. Let $\mf A$, $\mf B_t$, $t\ge 0$, stand for the
operators $\gamma^{-1} \Delta$, $\sqrt{2 \gamma^{-1} \chi(\rho(t,u))}
\nabla$, respectively, where $\rho$ is the solution of the heat
equation \eqref{f8} and $\chi$ is the compressibility given by
$\chi(\alpha) = \alpha (1-\alpha)$.

According to Holley-Stroock \cite{HS} theory of generalized
Ornstein-Uhlenbeck processes and to Stroock and Varadhan \cite{SV},
there exists a unique process $Z_t$ in $C([0,+\infty), \mathcal
S'(\mathbb R))$ with the following two properties:
$Z_0$ is a centered Gaussian field with covariance given by
\begin{equation}
\label{f13}
E[Z_0(G)Z_0(H)] = \int_{\bb R} G(u) H(u)  \chi(\rho_0(u))\, du
\end{equation}
for all $G$, $H$ in $\mc S(\bb R)$. Moreover, the processes $M_t(G)$,
$m_t(G)$ defined by
\begin{equation}
\label{f12}
Z_t(G) - Z_0(G) - \int_0^t Z_s(\mf A G) \, ds
\quad\text{and}\quad (M_t(G))^2 - \int_0^t 
\Vert \mf B_s G\Vert^2 \, ds
\end{equation}
are martingales with respect to the canonical filtration $\{ \mc F_s :
s\ge 0\}$ for all $G$ in $\mathcal S(\mathbb R)$. Of course, it is
enough to check these conditions for a dense family of functions in
$\mathcal S(\mathbb R)$.

Recall the definition of the class $\mathcal S_0(\mathbb R)$ and fix a
function $G$ in $\mc S_0(\bb R)$.  An elementary computation of the
characteristic function $\mathbb E_{\nu_{\rho_0(\cdot)}^N} [\exp\{i
\theta Z_0^N(G)\}]$ shows that $Z_0^N$ converges to a centered
Gaussian field with covariances given by \eqref{f13}.

Recall from \eqref{f14} the definition of the martingale $M^N_t(G)$
and fix a bounded function $U$ in $\mc F_s$.  To prove that $M_t(G)$
is a martingale, it is enough to show that
\begin{equation}
\label{f18}
\lim_{N\to\infty} \bb E_{\nu_{\rho_0(\cdot)}^N} [M^N_t(G) U]
\;=\; E [M_t(G) U] 
\end{equation}
for all $t\ge s$.

Since $Z_t^N(G)$ is bounded in $L^2 (\bb P_{\nu_{\rho_0(\cdot)}^N})$,
\eqref{f18} holds with $Z_t^N(G) - Z_0^N(G)$, $Z_t(G) - Z_0(G)$ in
place of $M^N_t(G)$, $M_t(G)$. By Schwarz inequality and a previous
estimate,
\begin{equation*}
\mathbb E_{\nu_{\rho_0(\cdot)}^N} \Big[ \Big(\frac{T_{\xi,G}}
{ T_{\xi,g}} \int_0^t Y_s^N (\Delta_N g_l) \, ds\Big)^2\Big] \;\le\;
\frac{C t^{2} N}{l_N^2} \Big(\frac{T_{\xi,G}} { T_{\xi,g}} \Big)^2
\end{equation*}
vanishes as $N\uparrow\infty$ for all $\xi$ in $\Omega_0$. On the
other hand, by Schwarz inequality, Theorem \ref{s11} and Lemma
\ref{s3},
\begin{eqnarray*}
\!\!\!\!\!\!\!\!\!\!\!\!\! &&
\mathbb E_{\nu_{\rho_0(\cdot)}^N} \Big[ \Big( \int_0^t 
\{ Y_s^N (\Delta_N G) - Z_s^N(\Delta_N G) \} \, ds
\Big)^2\Big] \\
\!\!\!\!\!\!\!\!\!\!\!\!\! && \quad
\le\; C_1 t^{5/2} \Big( \frac{1}{N} \sum_{x \in \mathbb Z}
|\Delta_N  G(x/N) - \gamma^{-1} (T_{\xi,l} \Delta_N G)(x)|\Big)^2\\
\!\!\!\!\!\!\!\!\!\!\!\!\! && \qquad
+ \; t^{2} \frac{1}{N} \sum_{x \in \mathbb Z}
\big\{\Delta_N  G(x/N) - \gamma^{-1} (T_{\xi,l} \Delta_N G)(x)\big\}^2
\end{eqnarray*}
vanishes as $N\uparrow\infty$ for all $\xi$ in $\Omega_0$. Replacing
$\Delta_N G$ by $\Delta G$ and recalling all previous estimates, we
deduce that
\begin{eqnarray*}
\lim_{N\to\infty} \bb E_{\nu_{\rho_0(\cdot)}^N} \Big [U \int_0^t
\gamma_1^N(s,G) \, ds \Big ] &\!\!\!=\!\!\!& \lim_{N\to\infty} \bb
E_{\nu_{\rho_0(\cdot)}^N} \Big [ U \int_0^t  \gamma^{-1} Z_s^N(\Delta
G) \, ds \Big ] \\
&\!\!\!=\!\!\!&  E \Big [ U \int_0^t  Z_s (\mf A  G) \, ds \Big ]
\end{eqnarray*}
because $Z_s^N (\Delta G)$ is bounded in $L^2$. This concludes the
proof of \eqref{f18}.

To prove \eqref{f18} with $m_t(G)$, $m^N_t(G) = M^N_t(G)^2 - \<M^N
(G)\>_t$ in place of $M_t(G)$, $M^N_t(G)$, observe first that $\bb
E_{\nu_{\rho_0(\cdot)}^N} [M^N_t(G)^4]$ is finite in view of Theorem
\ref{s11}, so that $\bb E_{\nu_{\rho_0(\cdot)}^N} [U M^N_t(G)^2]$
converges to $E [U M_t(G)^2]$.

To show that $\bb E_{\nu_{\rho_0(\cdot)}^N} [U \<M^N(G)\>_t]$
converges to $E[U \int_0^t \Vert \mf B_s G\Vert^2 ds]$, notice that
$\bb E_{\nu_{\rho_0(\cdot)}^N} [\<M^N(G)\>_t^2]$ is bounded uniformly
in $N$, for all $\xi$ in $\Omega_0$, and that $\<M^N(G)\>_t$ can be
written as
\begin{equation*}
\int_0^t ds\, \frac{1}{\gamma^2 N} \sum_{x \in \mathbb Z}
\xi_x^{-1} (\nabla_N G) (x)^2(\eta_s^N(x+1)-\eta_s^N(x))^2 
\end{equation*}
plus a remainder which vanishes as $N\uparrow\infty$ for all $\xi$ in
$\Omega_0$. By Theorem \ref{s11} and Schwarz inequality,
\begin{equation*}
\int_0^t ds\, \frac{1}{N} \sum_{x \in \mathbb Z}
\xi_x^{-1} (\nabla_N G) (x)^2 \{ \eta_s^N(x) - \rho_s^N(x)\}
\end{equation*}
vanishes in $L^2$, as well as the same expression with $\bar
\eta_s^N(x) \bar\eta_s^N(x+1)$ in place of $\bar \eta_s^N(x) =
\eta_s^N(x) - \rho_s^N(x)$. The penultimate integral is thus equal to
\begin{equation*}
\int_0^t ds\, \frac{1}{\gamma^2 N} \sum_{x \in \mathbb Z}
\xi_x^{-1} (\nabla_N G) (x)^2 \{\rho_s^N(x+1) + \rho_s^N(x) -
2 \rho_s^N(x+1) \rho_s^N(x) \}
\end{equation*}
plus a remainder which vanishes in $L^2 (\bb P_{\nu_{\rho_0(\cdot)}^N}
)$. This shows that $\bb E_{\nu_{\rho_0(\cdot)}^N} [U \<M^N(G)\>_t]$
converges to $E[U \int_0^t \Vert \mf B_s G\Vert^2 ds]$ and concludes
the proof of uniqueness.

\smallskip
\noindent{\bf Proof of Theorem \ref{t2}.}
By Lemma \ref{s3} and Theorem \ref{s11}, for each $t\ge 0$ and $G$ in
$\mc S(\bb R)$, $Z^N_t(G) - Y^N_t(G)$ vanishes in $L^2$ $\xi$-almost
surely as $N\uparrow\infty$. In particular, we may deduce from the
central limit theorem for $Z^N_t$ the convergence of the finite
dimensional distributions of $Y^N_t$. \qed \medskip

We conclude this section with the

\noindent{\bf Proof of Lemma \ref{s3}.}
Fix $G$ in $\mathcal S(\bb R)$ and recall that $\hat \xi^{-1}_x =
\xi_x^{-1} - \gamma$. For $N$ fixed, $T_\xi G(x)$ is well defined
because $\sum_{-k \le x\le k} \xi_x^{-1} (\nabla_N G)(x/N)$ is a
Cauchy sequence in $L^2(P)$.

By Doob's inequality, for each $x<y$, $A>0$,
\begin{equation*}
P \Big[ \max_{x < z\le y}  \Big\vert \frac 1N \sum_{w=x+1}^z 
\hat \xi^{-1}_w (\nabla_N G) (w/N) \Big\vert > A \Big] 
\;\leq\; \frac {C_0(G) E [\xi_x^{-6}]}{A^6 N^3} 
\end{equation*}
for some finite constant $C_0$ depending on $G$. Take $A=N^{-\{(1/4) +
  \epsilon\}}$ for some $0 < \epsilon < 1/12$, estimate the right hand
side by $C_0(G) E [\xi_x^{-6}] N^{6 \epsilon - 3/2}$ and let
$y\uparrow\infty$ to conclude by Borel-Cantelli that
\begin{equation*}
N^{1/4} \max_{x < z}  \Big\vert \frac 1N \sum_{w=x+1}^z 
\hat \xi^{-1}_w (\nabla_N G) (w/N) \Big\vert
\end{equation*}
vanishes, as $N\uparrow\infty$, almost surely. This proves the first
statement of the lemma.

Since 
\begin{equation*}
T_{\xi,G} \;=\; \frac 1N \sum_{x\in\bb Z} \xi_x^{-1} (\nabla_N G)(x/N)
\;=\; \frac 1N \sum_{x\in\bb Z} \hat \xi_x^{-1} (\nabla_N G)(x/N)\; ,
\end{equation*}
$T_{\xi,G}$ is absolutely bounded by $\sup_{x \in \bb Z} |T_\xi G(x)
-\gamma G(x/N)|$ and the second claim of the lemma follows from the
first one.

To prove the last two statement, notice that $(T_{\xi,l}G)(x) - \gamma
G(x/N)$ is absolutely bounded by $R_{\xi,G}^N (x) + |T_{\xi,G}| \, \mb
1\{0<x<l_N N\}$, where $R_{\xi,G}^N (x)$ is equal to
\begin{equation*}
\frac{1}{N} \Big\vert \sum_{y < x} \hat \xi^{-1}_y \nabla_N G(y/N) 
\Big\vert  \quad \text{for $x\le 0$} \;, \quad
\frac{1}{N} \Big\vert \sum_{y \geq x} \hat \xi^{-1}_x \nabla_N G(x/N)
\Big\vert \quad \text{for $x>0 $}\; .
\end{equation*}
In particular, by the first part of the proof and since $l_N =
N^{1/4}$, $\sup_{x \in \bb Z} |T_{\xi,l} G(x) -\gamma G(x/N)|$, $l_N
T_{\xi,G}$ vanishes, as $N\uparrow\infty$, $\xi$-almost surely.
On the other hand, by Tchebychev and H\"older inequality,
\begin{equation*}
P \Big[ \frac{1}{N} \sum_{x \in \bb Z} R_{\xi,G}^N (x)  > A \Big] 
\;\le\; \frac {C_0} {A^4 N^4} \sum_{x\in\bb Z} (1+|x|^{3(1+\epsilon)})
\, E\Big[ R_{\xi,G}^N (x) ^4 \Big]
\end{equation*}
for some $\epsilon >0$ and some finite constant $C_0 = C_0(\epsilon)$.
Since $G$ belongs to $\mc S(\bb R)$, the previous expectation is less
than or equal to $C_0 E[\xi_0^{-4}] N^{-2} F_G(x/N)$ for some rapidly
decreasing positive function $F_G$. The left hand side is thus bounded
above by $C_0(\epsilon) E[\xi_0^{-4}] N^{3\epsilon -2} A^{-4}$.
Choosing $0< \epsilon < 1/7$, $A=N^{-\epsilon}$ we conclude the proof
of the last statement of the lemma with a Borel-Cantelli argument.
\qed

\section{Central limit theorem for a tagged particle.}
\label{sec5}

We prove in this section Theorem \ref{t3}. Unless otherwise stated, we
assume throughout this section that $\rho_0$ is an initial condition
with first derivative in $L^1(\bb R) \cap L^\infty(\bb R)$ and second
derivative in $L^\infty(\bb R)$, and that the environment satisfies the
assumptions of the previous section. The proof follows closely the
approach presented in \cite{JL}. We omit therefore some details.

We first consider the current through a bond.  For each $x$ in
$\mathbb Z$, denote by $J_{x,x+1}^N (t)$ the current over the bond
$\{x, x+1\}$ in the time interval $[0,t]$. This is the total number of
particles which jumped from $x$ to $x+1$ minus the total number of
particles which jumped from $x+1$ to $x$ in the time interval $[0,t]$.

The current $J_{x,x+1}^N(t)$ can be related to the occupation
variables $\eta_t(x)$ through the formula
\begin{equation}
\label{f20}
J_{x-1,x}^N(t)-J_{x,x+1}^N(t) = \eta_t(x) - \eta_0(x)\; .
\end{equation} 
The first result states a law of large numbers for the current through
a bond assuming that the environment satisfies condition \eqref{f9}.

\begin{proposition}
\label{s5}
Consider a sequence $\{\xi_x : x\in \bb Z\}$ satisfying \eqref{f9} and
a profile $\rho_0 : \bb R_+\to [0,1]$ satisfying the assumptions
stated at the beginning of this section.  For every $\delta>0$, 
\begin{equation*}
\lim_{N\to\infty} \bb P_{\nu^N_{\rho_0(\cdot)}}^N \Big[ \, \Big\vert
\frac{J_{0,1}^N (t)}{N} + \int_0^t \gamma^{-1} (\partial_u \rho) (s,
0) \, ds \Big\vert > \delta \Big] \;=\; 0\;, 
\end{equation*}
where $\rho(t,u)$ is the solution of \eqref{f8}.
\end{proposition}

\begin{proof}
Fix $a >0$. Identity \eqref{f20} and a summation by parts give
that  
\begin{equation}
\label{f21}
\frac{1}{a N^2} \sum_{x =1}^{a N} \xi_x^{-1} 
\Big\{ J_{x,x+1}^N(t) - J_{0,1}^N(t) \Big\}
\;=\; \frac{1}{a N^2} \sum_{x =1}^{N a} 
\{\eta_t^N(x) -\eta_0^N(x)\} \sum_{k=x}^{a N} \xi_k^{-1} \;.
\end{equation}
Since the right hand side is of order $a$, the law of large numbers
for $J_0^N(t)/N$ follows from a law of large numbers for $a^{-1}
N^{-2} \sum_{x=1}^{Na} \xi_x^{-1} J_{x,x+1}^N(t)$.  We may rewrite
this latter expression as
\begin{equation}
\label{f29}
\frac{1}{a N^2} \sum_{x=1}^{a N} \xi_x^{-1} M_{x,x+1}^N(t)
\;+\; \frac{1}{a} \int_0^t \{ \eta_s^N(1) -\eta_s^N (a
N+1) \} \, ds \;,
\end{equation}
where 
\begin{equation*}
M_{x,x+1}^N(t) =: J_{x,x+1}^N(t) - N^2 \int_0^t
\xi_x \{ \eta_s^N(x)-\eta_s^N(x+1)\} \, ds\;,
\end{equation*}
$x$ in $\bb Z$, are orthogonal martingales with quadratic variation
$\< M_{x,x+1}^N\>_t$ given by
\begin{equation*}
N^2 \int_0^t \xi_x \{ \eta_s^N(x)-\eta_s^N(x+1)\}^2 \, ds\;.
\end{equation*}

In view of \eqref{f9} and of the explicit expression for the quadratic
variation of the orthogonal martingales $M_{x,x+1}^N(t)$, the first
term in \eqref{f29} vanishes in $L_2(\mathbb
P_{\nu^N_{\rho_0(\cdot)}}^N)$ as $N\uparrow\infty$. On the other hand,
by Theorem \ref{s11}, the variance of the second term in \eqref{f29}
vanishes as $N\uparrow\infty$. Its expectation is equal to
\begin{equation*}
\frac 1 a \int_0^t \{ \rho_s^{N,\xi} (1) - \rho_s^{N, \xi} 
(a N +1)\} \, ds \;.
\end{equation*}
By Lemma \ref{s7}, this integral converges to $a^{-1} \int_0^t \{
\rho_s (0) - \rho_s (a)\} \, ds$, where $\rho$ is the solution of
\eqref{f8}.  It remains to let $a \downarrow 0$ to conclude the proof.
\end{proof}

We prove now a quenched nonequilibrium central limit theorem for the
current. Let $\bar J_{x,x+1}^N(t) = J_{x,x+1}^N(t) -
E_{\nu^N_{\rho_0(\cdot)}} [ J_{x,x+1}^N(t) ]$.

\begin{proposition}
\label{s8}
There exists a total measure set $\Omega_0 \subset\Omega$ with the
following property. For each $\xi$ in $\Omega_0$, each $k \ge 1$ and
each $0\le t_1 <\cdots <t_k$, the random vector $N^{-1/2} (\bar
J_{-1,0}^N(t_1) , \dots, \bar J_{-1,0}^N(t_k))$ converges in law to a
Gaussian vector $(J_{t_1}, \dots, J_{t_k})$ with covariances given by
\begin{eqnarray*}
E[J_s J_t] &=&
\int_{-\infty}^0 dv \,  P[B_s\le v] \, P[B_t\le v] \, \chi(\rho_0(v)) \\
&+& \int_0^{\infty} dv \,  P[B_s\ge v] \, P[B_t\ge v] \, \chi(\rho_0(v)) \\
&+& 2 \gamma^{-1} \int_0^s dr \int_{-\infty}^{\infty} dv \, p_{t-r}(0,v) 
\, p_{s-r}(0,v) \, \chi(\rho(r,v))
\end{eqnarray*}
provided $s\le t$. In this formula, $B_t = B^0_{t/\gamma}$, where
$B^0_t$ is a standard Brownian motion starting from the origin, and
$p_t(v,w)$ is the kernel of $B_t$.
\end{proposition}

\begin{proof}
The proof of this proposition is similar to the one of Theorem 2.3 in
\cite{JL}. Some details are therefore omitted.

Let $H_0(u) = \mb 1\{u\ge 0\}$ and define the sequence $\{G_n : n\ge
1\}$ of approximations of $H_0$ by
\begin{equation*}
G_n(u) = \{1 - (u/n)\}^+ \mb 1\{u\ge 0\}\;.
\end{equation*}
We claim that for every $t\ge 0$,
\begin{equation}
\label{f30}
\lim_{n \to \infty} \bb E_{\nu^N_{\rho_0(\cdot)}} 
\Big[ N^{-1/2} \overline{J}_{-1,0}(t)
- Y^N_t(G_n) + Y^N_0(G_n) \Big]^2 \;=\; 0
\end{equation}
uniformly in N. The proof of \eqref{f30} relies on the estimates of
the two point space-time correlation functions, presented in Lemma
\ref{s13}, and follows closely the proof of Proposition 3.1 in
\cite{JL}. We leave the details to the reader.

Fix $t\ge 0$ and $n\ge 1$. By approximating $G_n$ in $L^2(\bb R) \cap
L^1(\bb R)$ by a sequence $\{H_{n,k} : k\ge 1\}$ of smooth functions
with compact support, recalling Theorem \ref{t2}, we show that
$Y^N_t(G_n)$ converges in law to a Gaussian variable denoted by
$Y_t(G_n)$.

By \eqref{f30}, $\{Y^N_t(G_n) - Y^N_0(G_n) : n\ge 1\}$ is a Cauchy
sequence uniformly in $N$. In particular, $Y_t(G_n) - Y_0(G_n)$ is a
Cauchy sequence and converges to a Gaussian limit denoted by $Y_t
(H_0) - Y_0(H_0)$. Therefore, by \eqref{f30}, $N^{-1/2}
\overline{J}_{-1,0}(t)$ converges in law to $Y_t (H_0) - Y_0(H_0)$.

The same argument show that any vector $N^{-1/2}
(\overline{J}_{-1,0}(t_1), \dots, \overline{J}_{-1,0}(t_k))$ converges
in law to $(Y_{t_1} (H_0) - Y_0(H_0), \dots, Y_{t_k} (H_0) -
Y_0(H_0))$. The covariances can be computed since by (\ref{f17})
\begin{eqnarray*}
\!\!\!\!\!\!\!\!\!\!\!\!\!\!\! &&
E\Big[ \big\{ Y_{t} (H_0) - Y_0(H_0) \big\} \big\{ Y_{s} (H_0) -
Y_0(H_0) \big\}  \Big] \\
\!\!\!\!\!\!\!\!\!\!\!\!\!\!\! && \quad
=\;\lim_{n\to\infty} E\Big[ \big\{ Y_{t} (G_n) - Y_0(G_n) \big\}
\big\{ Y_{s} (G_n) - Y_0(G_n) \big\}  \Big] \\
\!\!\!\!\!\!\!\!\!\!\!\!\!\!\! && \quad
=\; \lim_{n\to\infty} \Big\{ \int_{\bb R} \big\{ (T_{t} G_n) (T_s G_n)  +
G_n^2  - (T_t G_n) G_n  - (T_s G_n) G_n \big\} \chi(\rho_0 (u)) \\
\!\!\!\!\!\!\!\!\!\!\!\!\!\!\! && \qquad\qquad\qquad
\;+\; 2 \gamma^{-1} \int_0^s dr \int_{\bb R} (\nabla T_{t-r} G_n) \,
(\nabla  T_{s-r} G_n) \chi(\rho (r,u)) \Big\}\;.
\end{eqnarray*}
A long but elementary computation permits to recover the expression
presented in the statement of the proposition.
\end{proof}

We turn now to the behavior of a tagged particle.  Let
$\nu_{\rho_0(\cdot)}^{N,\ast}$ be the product measure
$\nu_{\rho_0(\cdot)}^N$ conditioned to have a particle at the origin.
All our previous results for the process starting from
$\nu_{\rho_0(\cdot)}^N$ remain in force for the process starting from
$\nu_{\rho_0(\cdot)}^{N,\ast}$, since we can couple both processes in
such a way that they differ at most at one site at any given time.

Denote by $X_t^N$ the position at time $t \geq 0$ of the particle
initially at the origin.  Since the relative ordering of particles is
conserved by the dynamics, a law of large numbers for $X_t^N$ is a
consequence of the hydrodynamic limit and the law of large numbers for
the current (\cite{RV}, \cite{LV}, \cite{JL}).  In fact, the
distribution of $X_t^N$ can be obtained from the joint distribution of
the current and the empirical measure via the relation
\begin{equation}
\label{f27}
\{ X_t^N \geq n \} \;=\;
\{ J_{-1,0} (t) \geq \sum_{x=0}^{n-1} \eta_t(x) \} 
\end{equation}
for all $n\ge 0$ and a similar relation for $n\le 0$.  

\begin{theorem}
\label{s14}
Consider a sequence $\{\xi_x : x\in \bb Z\}$ satisfying \eqref{f9} and
a profile $\rho_0 : \bb R_+\to [0,1]$ satisfying the assumptions
presented at the beginning of this section. For every $t\ge 0$,
$X_t/N$ converges in $\bb
P_{\nu_{\rho_0(\cdot)}^{N,\ast}}$-probability to $u_t$, the solution
of
\begin{equation}
\label{f33}
\int_0^{u_t} \rho(t,u) du = - \frac 1 \gamma \int_0^t
(\partial_u\rho)(s,0) ds\;.
\end{equation}
\end{theorem}

Notice that $u_t$ satisfies the differential equation
\begin{equation*}
\dot u_t \;=\; - \frac 1\gamma
\frac{(\partial_u \rho)(t,u_t)}{\rho(t,u_t)}\;\cdot
\end{equation*}
The proof of this result is similar to the one of Theorem 2.5 in
\cite{JL} and left to the reader. 

It remains to prove a central limit theorem for the position of the
tagged particle.

\noindent{\bf Proof of Theorem \ref{t3}.}
Recall the definition of $u_t^N$ presented just before the statement
of the theorem, assume that $u^N_t> 0$ and fix $a$ in $\bb R$.
By equation (\ref{f27}), the set $\{X_t \ge u^N_t + a\sqrt{N}\}$ is
equal to the set in which
\begin{equation}
\label{eq:03}
\overline{J}_{-1,0} (t) \geq
\sum_{x=0}^{u^N_t} \overline{\eta}_t (x)
\;+\; \sum_{x=1}^{a\sqrt{N}  -1} \eta_t(x+u^N_t)
\;-\;  \Big\{ \mathbb{E}_{\nu^{N,*}_{\rho_0(\cdot)}} [J_{-1,0}(t)]
-\sum_{x=0}^{u^N_t} \rho^{N,*}_t(x)]\Big\} \;.
\end{equation}
We claim that second term on the right hand side of this equation
divided by $\sqrt{N}$ converges to its mean in $L^2$. Indeed, by
Theorem \ref{s11}, its variance is bounded by $C_0(\eps, \rho_0) a
N^{-1/2}$ for some finite constant $C_0$. Notice that we are taking
expectations with respect to a measure, $\nu^{N,*}_{\rho_0(\cdot)}$,
whose associated profile does not have a bounded first derivative.
However, coupling this measure with $\nu^{N}_{\rho_0(\cdot)}$, in such
a way that they differ at most by one particle at every time, we can
still show that the variance is bounded by $C_0(\eps, \rho_0) a
N^{-1/2}$ as claimed. The same ideas, the linearity of equation
\eqref{f16} and Nash estimate, stated in Proposition \ref{s9} below,
show that $\rho^{N,*}$ converges uniformly on compact sets to the
the solution of the heat equation \eqref{f8} because $\rho^{N}$
converges in view of Lemma \ref{s7}. 

To compute the expectation of the second term on the right hand side
of \eqref{eq:03}, observe that the middle term in \eqref{f32} is equal
to $\mathbb{E}_{\nu^{N,*}_{\rho_0(\cdot)}} [J_{-1,0}(t)]$. By the
proof of the law of large numbers for the current, this middle
expression divided by $N$ converges to $- \gamma^{-1} \int_0^t
(\partial_u \rho)(s,0) ds$. In particular, by the law of large numbers
for the empirical measure and by relation \eqref{f33}, $N^{-1} u^N_t $
converges to $u_t$. Hence, by the uniform convergence of $\rho^{N,*}$,
$$
\frac 1{\sqrt{N}} \sum_{x=1}^{a\sqrt{N}  -1} \rho^{N,*}_t(x+u^N_t)
$$
converges to $a \rho(t,u_t)$ and so does in probability the second
term on the right hand side of (\ref{eq:03}).

By definition of $u^N_t$, the third term on the right hand side is
absolutely bounded by $1$.

Finally, by \eqref{f30}, for fixed $t$, $N^{-1/2} \{
\overline{J}_{-1,0} (t) - \sum_{x=0}^{u^N_t} \overline{\eta}_t (x)\}$
behaves as $Y^N_t(G_n) - Y^N_0(G_n) - Y^N_t(\mb 1\{[0,u^N_t/N]\})$, as
$N\uparrow\infty$, $n\uparrow\infty$. Repeating the arguments
presented in the proof of Proposition \ref{s8}, since $u^N_t/N$
converges to $u_t$, we show that this latter variable converges in law
to a centered Gaussian variable, denoted by $W_t$, and which is
formally equal to $Y_t(H_{u_t})-Y_0(H_0)$, where $H_{a}(u) = \mb
1\{u\ge a\}$.

Up to this point we proved that
$$
\lim_{N\to\infty} \bb P_{\nu^{N,*}_{\rho_0(\cdot)}} \Big[
\frac{X^N_t - u^N_t}{\sqrt{N}} \ge  a\Big] \;=\;
P[W_t \ge a \rho(t,u_t)]
$$
provided $u_t>0$. Analogous arguments permit to prove the same
statement in the case $u_t=0$, $a>0$. By symmetry around the origin,
we can recover the other cases: $u_t<0$ and $a$ in $\bb R$, $u_t=0$
and $a<0$.
                                                 
Putting all these facts together, we conclude that for each fixed $t$,
$(X_t- u^N_t)/\sqrt{N}$ converges in distribution to the Gaussian
$W_t/\rho(t,u_t) = [Y_t(H_{u_t})-Y_0(H_0)]/\rho(t,u_t)$.  The same
arguments show that any vector $(N^{-1/2} [X_{t_1}- u^N_{t_1}], \dots,
N^{-1/2} [X_{t_k}-u^N_{t_k}])$ converges to the corresponding centered
Gaussian vector. It remains to compute the covariances, which can be
derived as in the proof of Proposition \ref{s8}. Details are left
to the reader. \qed

\section{Correlation estimates}
\label{sec4}

We assume throughout this section that $\{\xi_x : x\in\bb Z\}$ is a
sequence of numbers bounded below and above: $0 < \eps < \xi_x <
\eps^{-1}$ for all $x$, and that the profile $\rho_0: \bb R \to [0,1]$
has bounded first derivative. Recall that $\rho_t^N(x) = \bb
E_{\nu_{\rho_0^N(\cdot)}} [\eta_t(x)]$ satisfies equation
\eqref{f16}.

For $n \geq 1$, denote by $\mc E_n$ the subsets of $\bb Z$ with $n$
points. For each $\x_n=\{x_1,...,x_n\}$ in $\mc E_n$, let 
\begin{equation*}
\varphi_t(\x_n) = \bb E_{\nu_{\rho_0^N(\cdot)}} \big[ \prod_{i=1}^n \{
\eta_t (x_i) -\rho_t^N(x_i)\}\big]\; .
\end{equation*}

\begin{theorem}
\label{s11}
Fix a finite time interval $[0,T]$ and an initial profile $\rho_0$
with bounded first derivative. There are constants $C_n$, depending
only on $\eps$, $\rho_0$, $n$ and $T$, such that
\begin{equation*}
\sup_{\substack{\x_{2n} \in \mc E_{2n} \\ t \in [0,T]}} |\varphi_t(\x_{2n})|
\leq \frac{C_{2n}}{N^n}, \qquad 
\sup_{\substack{\x_{2n+1} \in \mc E_{2n+1} \\ t \in [0,T]}} 
|\varphi_t (\x_{2n+1})| \leq \frac{C_{2n+1}\log N}{N^{n+1}}\; .
\end{equation*}
\end{theorem}

The proof of this theorem follows closely the proof of \cite{FPSV}
for the simple exclusion process without environment. We start with a
Nash estimate for the transition probability of a random walk in
elliptic environment \cite{K}, \cite{CKS}, \cite{Dav}. Denote by $\mc
L_1$ the generator of a random walk in the bond environment $\xi$:
\begin{equation*}
(\mc L_1 f)(x) = \xi_{x-1} \{f(x-1) - f(x)\} \;+\; \xi_{x}
\{f(x+1) - f(x)\}\;.
\end{equation*}
Let $p^\xi_t(x,y)$ be the transition probability associated to the
generator $\mc L_1$.

\begin{proposition}
\label{s9}
There exists a finite constant $C_0 (\epsilon)$, depending only on
$\epsilon$, such that $p^\xi_t(x,y) \leq C_0(\epsilon) t^{-1/2}$ for
all $x$, $y$ in $\bb Z$, $t\ge 0$.
\end{proposition}

The proof of Theorem \ref{s11} relies also on a comparison between
the semigroup associated to the evolution of $n$ exclusion particles
with the semigroup associated to $n$ independent particles.  For $n\ge
1$, denote by $\mc L_n$ the generator corresponding to the evolution
of $n$ exclusion particles in the environment $\xi$:
\begin{eqnarray*}
(\mc L_n h)(\x_n ) &=&   N^2 \sum_{i=1}^n
{\bf 1} \{\x+e_i \in \mc E_n\} \, \xi_{x_i} \, [h (\x_n +e_i) - h (\x_n)] \\
&+& N^2 \sum_{i=1}^n {\bf 1} \{\x-e_i \in \mc E_n\} \, \xi_{x_i-1}\,
[h(\x_n-e_i) - h(\x_n)] 
\end{eqnarray*}
for every function $h: \mc E_n\to\bb R$. In this formula, for $1\le
i\le n$, $e_i$ stands for the $i$-th canonical vector in $\bb R^n$ and
$\x_n$ is understood as the vector $(x_1, \dots, x_n)$.  Denote by
$S_n(t)$ the semigroup associated to $\mc L_n$ and by $S_n^0(t)$ the
semigroup associated to $n$ independent particles evolving in the
environment $\xi$.

A bounded symmetric function $f : \mathbb Z^2 \to \mathbb R$ is said
to be definite positive provided
\begin{equation*}
\sum_{x,y} f(x,y) \beta_x \beta_y \geq 0
\end{equation*}
for every sequence $\{ \beta_x : x\in \bb Z\}$ in $\ell_1(\bb Z)$ with
$\sum_x \beta_x=0$. A bounded symmetric function $f : \mathbb Z^n \to
\mathbb R$ is said to be definite positive if it is so for each pair
of coordinates.  From \cite{Lig} we have that

\begin{proposition}
\label{s10} 
Let $f : \mathbb Z^n \to \mathbb R$ be a bounded, symmetric, definite
positive function. Then, 
\begin{equation*}
S_n(t) f \leq S_n^0 (t) f
\end{equation*}
for all $t\ge 0$.
\end{proposition}

Theorem \ref{s11} is based on an induction argument.  Observe first
that
\begin{equation}
\label{f23}
\frac{d}{dt} \varphi_t(\x_n) \;=\; (\mc L_n \varphi_t) (\x_n )
\;+\; \Gamma_t (\x_n) \;,
\end{equation}
where
\begin{eqnarray*}
\Gamma_t (\x_n) &=&  2 N^2 \sum_{\substack{x\in \bb Z \\ x, x+1 \in \x_n}}  
\xi_{x}\, [\varphi_t (\x_n^{x+1}) - \varphi_t (\x_n^{x})]
\, [\rho_t^N(x+1)- \rho_t^N(x)] \\ 
&-& N^2 \sum_{\substack{x\in \bb Z \\ x, x+1 \in \x_n}}  \xi_{x}
\, \varphi_t(\x_n^{x,x+1})\, [\rho_t^N(x+1)-\rho_t^N(x)]^2 \;.
\end{eqnarray*}
Here and below $\x_n^{y}$, $\x_n^{y,z}$ stand for the configuration 
$\x_n \setminus \{y\}$, $\x_n \setminus \{y,z\}$, respectively.

In view of the differential equation \eqref{f23}, we can represent
$\varphi_t(\x_n)$ as an expectation with respect to a random walk in
an environment $\xi$ with sources at the boundary $\partial \mc E_n =
\{\x_n \in \mc E_n ; \min_{i\not = j} |x_i-x_j|=1\}$: Denote by $\bb
E_{\x_n}$ (resp. $\bb E_{\x_n}^0$) the expectation with respect to
$n$ exclusion (resp. independent) particles starting at $\x_n$. Since
$\varphi_0(\x_n) = 0$, we have that
\begin{equation}
\label{f24}
\varphi_t(\x_n) = \int_0^t ds\, \bb E_{\x_n} 
\big[ \Gamma_s (\x_n(t-s)) \big]\;.
\end{equation}

Since $\varphi_t(x)=0$ for all $x$ in $\bb Z$, $t\ge 0$, to start the
induction argument, set $n=2$ and remark that the first term in the
definition of $\Gamma_t$ vanishes. On the other hand, by \eqref{f28}
below, the derivative $\nabla_N \rho_t^N$ is uniformly bounded.  Since
the environment is elliptic and $\varphi_t(\phi)=1$, $\Gamma_t(\x)$ is
absolutely bounded by $C_0(\eps, \rho_0) \mb 1\{ \x_2 \in \partial \mc
E_2\}$ for some finite constant $C_0$.

The function $f:\bb Z^2 \to \bb R$ defined by $f(x,y) = 2 \mb 1\{x=y\}
+ \mb 1\{|x-y|=1\}$ is bounded, symmetric and definite positive.
Therefore, by Proposition \ref{s10}, by the integral representation
\eqref{f24} of $\varphi_t$ and by the previous estimate of $\Gamma_t$,
\begin{equation*}
\varphi_t(x,y) \leq C(\eps,\rho_0)\int_0^t 
\bb E_{(x,y)}^0 [ f(\x_2(s))] \, ds \; .
\end{equation*}
It remains to apply Proposition \ref{s9} and to integrate in time,
keeping in mind that time is speeded up by $N^2$, to obtain that
\begin{equation*}
\sup_{x \neq y \in \bb Z} |\varphi_t (x,y)| \leq 
\frac{C(\eps, \rho_0) \sqrt{t}}{N} \leq 
\frac{C(\eps, \rho_0, T)}{N}
\end{equation*}
for all $0\le t\le T$.

To extend this estimate to $n\ge 3$, we need to exploit the
non-trivial cancellations in the first term of the definition of
$\Gamma_t$. For $n\ge 1$, let
\begin{eqnarray*}
\!\!\!\!\!\!\!\!\!\!\!\!\!\!\!\! &&
A_t^n \;=:\; \sup_{\x_n \in \mc E_n} |\varphi_t(\x_n)|\;, \\
\!\!\!\!\!\!\!\!\!\!\!\!\!\!\!\! && \quad
B_t^n \;=:\; \sup_{x\in\bb Z} \sup_{\substack{ \x_{n-1} \in \mc E_{n-1} \\ 
\x_{n-1} \not\ni x, x+1}}
| \varphi_t (\x_{n-1} \cup\{x\})-\varphi_t(\x_{n-1} \cup\{x+1\} )|\; .
\end{eqnarray*}
We claim that there exists a finite sequence of constants $C(\eps,
\rho_0, n)$, $n\ge 2$, such that
\begin{eqnarray}
\label{f25}
\!\!\!\!\!\!\!\!\!\!\!\!\!\!\!\! &&
A_t^n \leq C(\eps, \rho_0, n) \int_0^t \big\{ 
N B_s^{n-1} + A_s^{n-2} \big\} \frac{ds}{N \sqrt{t-s}}\; , \\
\!\!\!\!\!\!\!\!\!\!\!\!\!\!\!\! && \quad
B_t^n \leq C(\eps, \rho_0, n) \int_0^t \big\{ N B_s^{n-1} 
+ A_s^{n-2}\big\} \frac{ds}{1+N^2(t-s)} \;.
\nonumber
\end{eqnarray}
Theorem \ref{s11} follows from these bounds and elementary
computations.

It remains to prove the estimates \eqref{f25}. The first one is
simpler and follows the same steps presented for $n=2$. Fix $n\ge 3$
and a configuration $\x_n$ in $\mc E_n$. Assume that the particles are
evolving according to a stirring process.  By \eqref{f24} and by
definition of $A_s^{k}$, $B_s^{k}$,
\begin{align*}
|\varphi_t (\x_n)| & \leq \int_0^t ds\, \bb E_{\x_n} 
\Big[ \, \big |\Gamma_s(\x_n (t-s)) \big | \, \Big] \\
& \leq C(\eps,\rho_0,n) \int_0^t ds\, \bb P_{\x_n}
\big [\x_n(t-s) \in \partial \mc E_n \big]
\{N B_s^{n-1} + A_s^{n-2}\} 
\end{align*}
for some finite constant $C(\eps,\rho_0,n)$. Since $\x_n(t-s)$ belongs
to the boundary of $\mc E_n$, there are at least two particles at
distance one. By definition of the the stirring process, any pair of
particles evolves according to a symmetric exclusion process in the
environment $\xi$. In particular, comparing the original process with
independent particles and applying Nash estimate, we can bound the
probability appearing in the last displayed formula by $C \{
N^2(t-s)\}^{-1/2}$.  This proves the first estimate in \eqref{f25}.

We now turn to $B_t^n$. Since $B_t^1 =0$, fix $n\ge 2$, $x$ in $\bb Z$
and $\x_{n-1}$ in $\mc E_{n-1}$ such that $x$, $x+1 \not\in \x_{n-1}$.
Consider $n+1$ particles evolving on $\bb Z$ according to the
following rules. They start from $\x_{n-1}$, $x$, $x+1$ and evolve
according to a stirring process. However, when the particles starting
at $x$ and $x+1$ are at distance $1$, each one jumps, independently
from the other, to the site occupied by the other at the rate
determined by the environment. Once these particles occupy the same
site, they remain together for ever.  Notice that the two
distinguished particles behave until they meet exactly as two
independent particles.

Denote by $\bb P_{\x_{n-1}, x, x+1}$, $\bb E_{\x_{n-1}, x, x+1}$ the
probability and the expectation corresponding to the evolution just
described.  Let $\tau$ be the coalescence time of the distinguished
particles and let $\x_n (t,x)$, $\x_n (t,x+1)$ be the configuration at
time $t$ of the system starting from $\x_{n-1} \cup \{x\}$, $\x_{n-1}
\cup \{x+1\}$, respectively. By construction, $\x_n (t,x) = \x_n
(t,x+1)$ for $t\ge \tau$. In particular,
\begin{eqnarray*}
\!\!\!\!\!\!\!\!\!\!\!\!\!\!\! &&
\varphi_t(\x_{n-1} \cup \{x\}) - \varphi_t(\x_{n-1} \cup \{x+1\}) \\
\!\!\!\!\!\!\!\!\!\!\!\!\!\!\! && \qquad
=\; \int_0^t ds\, \bb E_{\x_{n-1}, x, x+1} 
\big [ \Gamma_s (\x_n (t,x))-\Gamma_s (\x_n (t,x+1) ) \big] \\
\!\!\!\!\!\!\!\!\!\!\!\!\!\!\! && \qquad
=\; \int_0^t ds\, \bb E_{\x_{n-1}, x, x+1} 
\big [ \mb 1\{\tau > t-s\} \{ \Gamma_s (\x_n (t-s,x))
- \Gamma_s (\x_n (t-s,x+1) ) \} \big]\; .
\end{eqnarray*}
By definition of $A^k_s$ and $B^k_s$, this expression is less than or
equal to 
\begin{equation*}
C_0 \sum_{y=x,x+1} \int_0^t ds\, \{ N B_s^{n-1} + A_s^{n-2}\} 
\, \bb P_{\x_{n-1}, x, x+1} \Big[ \tau> t-s \,,\,
\x_n (t-s,y) \in \partial \mc E_n \Big]
\end{equation*}
for some finite constant $C_0 = C_0 (\eps,\rho_0,n)$.  In view of Nash
estimate, replacing the
indicator function ${\bf 1}\{ \tau > t-s\}$ by ${\bf 1} \{\tau >
(t-s)/2\}$ and applying the Markov property at time $(t-s)/2$, we
bound the previous expression by
\begin{equation*}
C_0 \int_0^t ds\, \{ N B_s^{n-1} + A_s^{n-2}\} 
\frac{1}{\sqrt{1+N^2 (t-s)}} \, \bb P_{\x_{n-1}, x, x+1} 
\big [ \tau> (t-s)/2  \big ]\; .
\end{equation*}
By \eqref{f26} below, the probability appearing in the previous
formula is bounded above by $C(\eps) \{ 1 + N^2(t-s)\}^{-1/2}$. This
concludes the proof of estimate \eqref{f25} and the one of Theorem
\ref{s11}. \qed

Let $x_t$ be a random walk in the environment $\{\xi_x : x\in \bb Z\}$
starting from $x_0=0$. Denote by $P$ the probability measure on the
path space $D(\bb R_+, \bb Z)$ induced by $x_t$. For each $a\not = 0$,
let $\tau_a$ be the first time the random walk $x_t$ reaches $a$:
\begin{equation*}
\tau_a =: \inf\{t \geq 0 ; x_t =a\}.
\end{equation*}

\begin{lemma}
\label{s12}
There exists a finite constant $C_0=C_0(\eps)$, depending only $\eps$,
such that 
\begin{equation*}
P(\tau_a > t) \leq \frac{C_0 a}{\sqrt {1+t}}
\end{equation*}
for all $t>0$.
\end{lemma}

\begin{proof}
Define the function $u: \bb Z \to \bb R$ by $u(0)=0$, $u(x+1)-u(x) =
\xi_x^{-1}$.  Since the environment is elliptic, $\eps \le u(x)/x \le
\eps^{-1} $ for all $x\not = 0$. Moreover, an elementary computation
shows that $u(x_t)$ is a martingale of quadratic variation
$\<u(x)\>_t$ given by 
$$
\int_0^t (\xi_{x_s-1}^{-1}+ \xi_{x_s}^{-1})\, ds \; .
$$

Fix $b<0<a$ and set $\tau = \min\{\tau_a, \tau_b\}$.  By Doob's
optional sampling theorem, $E[u(x_\tau)]=0$ and $E[u(x_\tau)^2 -
\<u(x)\>_\tau]=0$. Therefore,
\begin{equation*}
P(\tau_a< \tau_b) = \frac{-u(b)}{u(a)-u(b)}\; , \quad 
-u(a)u(b) = E \int_0^\tau (\xi_{x_s-1}^{-1}+ \xi_{x_s}^{-1})ds\; ,
\end{equation*}
so that $E[\tau] \leq -u(a)u(b)(2\eps)^{-1}$. In particular,
\begin{equation*}
P(\tau_a > t) \leq P(\tau >t ) + P(\tau_a > \tau_b) 
\leq \frac{-u(a)u(b)}{2\eps t} +  \frac{u(a)}{u(a)-u(b)}.
\end{equation*}
Minimizing over $b>0$ we conclude the proof of the lemma.
\end{proof}

The same ideas provide a bound on the coalescence time of two
independent particles in the environment $\xi$. Fix $x$ in $\bb Z$ and
consider two independent random walks $X_t$, $Y_t$, on the environment
$\xi$ such that $X_0=x$, $Y_0=x+1$. For $b>0$, let $\tau^*$, $\tau_b$
be the first time such that $Y_t=X_t$, $Y_t=X_t +b$, respectively.

Recall the definition of the function $u$ defined in the proof of
Lemma \ref{s12}. Since $X_t$, $Y_t$ are independent, $M_t = u(Y_t) -
u(X_t) - 1$ is a martingale. Repeating the arguments presented in the
proof of Lemma \ref{s12}, we obtain that
\begin{equation}
\label{f26}
P(\tau^* > t) \leq \frac{C_0}{\sqrt {1+t}}
\end{equation}
for all $t>0$ and some finite constant $C_0$ depending only on $\eps$.
Of course, when the time is speeded up by $N^2$, $t$ is replaced by
$tN^2$.  \medskip

A bound on the space-time correlations can be deduced from Theorem
\ref{s11}. For $x$, $y$ in $\bb Z$ and  $s\le t$, let
\begin{equation*}
\psi_{s,t}(y;x)  \;=\; \bb E_{\nu_{\rho_0^N(\cdot)}} 
\big[ \{ \eta_s(y) - \rho_s^N(y)\}
\{\eta_t(x) - \rho_t^N(x)\} \big] \; .  
\end{equation*}

\begin{lemma}
\label{s13}
There exists a finite constant $C_0$, depending only on $\eps$,
$\rho_0$ such that 
\begin{equation*}
\sup_{x,y \in \bb Z} |\psi_{s,t}(y;x)| \;\leq\; \frac{C_0}{N} 
\Big\{ \sqrt s + \frac{1}{\sqrt{t-s}} \Big\}
\end{equation*}
for all $0\le s\le t$.
\end{lemma}

\begin{proof}
Fix $s\ge 0$ and $y$ in $\bb Z$. For $t\ge s$, $x$ in $\bb Z$, let
$\psi_t (x) = \psi_{s,t}(y;x)$.  Notice that $\psi_t$ satisfies the
Cauchy problem 
$$
\begin{cases}
\frac{d}{dt} \psi_{t}(x) & = \mc L_1 \psi_{t}(x) \\
\psi_s(x) &= {\bf 1} \{x \neq y\} \varphi_s (x,y) 
+ {\bf 1}\{x=y\} \rho_s^N(y) (1-\rho_s^N(y)) \;, 
\end{cases}
$$
where $\mc L_1$ is the generator defined at the beginning of this
section. It remains to recall Nash estimate for the semigroup and the
proof of Theorem \ref{s11}, in which we showed that $\varphi_s (x,y)$
is bounded by $C \sqrt{s}/N$.
\end{proof}

We conclude this section with a result on the solution of the discrete
linear equation \eqref{f16}.

\begin{lemma}
\label{s7}
Let $\rho_0: \bb R \to [0,1]$ be a profile whose first derivative
$\rho_0'$ belongs to $L^1(\bb R) \cap L^\infty(\bb R)$ and whose
second derivative $\rho_0''$ belongs to $L^\infty(\bb R)$.  The
solution $\rho_t^N$ of equation \eqref{f16} converges uniformly on
compact sets of $\bb R_+\times \bb R$ to the solution of \eqref{f8}.
In particular, for all $t \ge 0$ and all function $G$ in
$C^{1}_{0}(\bb R)$,
\begin{eqnarray*}
\!\!\!\!\!\!\!\!\!\!\!\!\!\!\! &&
\lim_{N \to \infty} 
\frac{1}{N} \sum_{x \in \mathbb Z} G (x/N) \rho_t^N(x) 
\;=\; \int G (u) \rho (t,u) du \;, \\
\!\!\!\!\!\!\!\!\!\!\!\!\!\!\! &&
\lim_{N \to \infty} \frac{1}{N} \sum_{x \in \mathbb Z} 
G (x/N) (\nabla_N \rho_t^N) (x)  \;=\; 
\int G (u) (\partial_u \rho) (t,u) du \;,
\end{eqnarray*}
where $\rho_t(u)$ is the solution of linear heat equation \eqref{f8}. 
\end{lemma}

\begin{proof}
Consider the initial condition $\rho^\xi_0 = \rho^{N, \xi}_0 : \bb Z
\to \bb R$ defined by $\rho^\xi_0 (0) = \rho^N_0 (0)$, $(\nabla_N
\rho^\xi_0) (x) = \xi_x^{-1} (\nabla_N \rho^N_0) (x)$.  By the
estimates presented at the beginning of Section \ref{sec3},
$\rho_0^\xi (x) - \gamma \rho_0 (x)$ vanishes, as $N\uparrow\infty$,
uniformly in $x$.

Denote by $\rho_t^\xi = \rho_t^{N,\xi}$ the solution of equation
\eqref{f16} with initial condition $\rho_0^\xi$. We claim that the
sequence $\{\rho_t^{N,\xi} : N\ge 1\}$ is equicontinuous on each
compact set of $\bb R_+\times \bb R$. The proof relies on uniform
bounds of $\rho_t^\xi$, $\nabla_N \rho_t^\xi$, $(d/dt) \rho_t^\xi$.

First of all,  by the maximum principle,
\begin{equation*}
\inf_{x \in \mathbb Z} \rho_0^\xi (x) \;\leq\; 
\inf_{x \in  \mathbb Z} \rho_t^\xi (x) \;\leq\; 
\sup_{x \in \mathbb Z} \rho_t^\xi (x) \;\leq\; 
\sup_{x \in \mathbb Z} \rho_0^\xi (x) \;.
\end{equation*}
Denote by $\nabla_\xi$ the discrete derivative defined by $(\nabla_\xi
h) (x) = N \xi_x \{h(x+1) - h(x)\}$.  $(\nabla_\xi \rho_t^\xi)$ satisfies
the equation
\begin{equation*}
\frac{d}{dt} (\nabla_\xi \rho_t^\xi) (x) \;=\; N^2 \xi_x \big \{
(\nabla_\xi \rho_t^\xi) (x+1) + (\nabla_\xi \rho_t^\xi) (x-1) -
2 (\nabla_\xi \rho_t^\xi) (x) \big \}\;.
\end{equation*}
In particular, $\nabla_\xi \rho_t^\xi$ satisfies the maximum principle
and is uniformly bounded because we assumed the initial condition to
have a bounded derivative.

Let $\dot \rho_t^\xi = (d/dt) \rho_t^\xi$. By definition, 
\begin{equation*}
\dot \rho_t^\xi(x) \;=\; N\{ (\nabla_\xi \rho_t^\xi) (x) 
- (\nabla_\xi \rho_t^\xi) (x-1) \} \;=\; 
(\nabla_N \nabla_\xi \rho_t^\xi) (x-1)\;.
\end{equation*}
Since
\begin{eqnarray*}
\!\!\!\!\!\!\!\!\!\!\!\!\! &&
\frac{d}{dt} \dot \rho_t^\xi(x) \;=\; 
N^2  \Big \{ \xi_x \Big[ (\nabla_N \nabla_\xi \rho_t^\xi) (x)
- (\nabla_N \nabla_\xi \rho_t^\xi) (x-1) \Big ] \\
\!\!\!\!\!\!\!\!\!\!\!\!\! && \qquad\qquad\qquad\qquad\qquad
- \; \xi_{x-1} \Big[ (\nabla_N \nabla_\xi \rho_t^\xi) (x-1) 
- (\nabla_N \nabla_\xi \rho_t^\xi) (x-2) \Big ] \Big\}\;,
\end{eqnarray*}
$\dot \rho_t^\xi(x) = (\nabla_N \nabla_\xi \rho_t^\xi) (x-1)$
satisfies a maximum principle. By definition of $\rho_0^\xi$,
$\nabla_\xi$,
\begin{equation*}
(\nabla_N \nabla_\xi \rho_0^\xi) (x-1) \;=\;
N \big\{ (\nabla_\xi \rho_0^\xi) (x)
- (\nabla_\xi \rho_0^\xi) (x-1) \big \} \;=\; (\Delta_N \rho_0) (x)
\end{equation*}
In particular, $(d/dt) \rho_t^\xi(x)$ is uniformly bounded because we
assumed the initial condition to have bounded second derivative.

Notice that the previous bound does not hold for the initial condition
$\rho_0$ since $\nabla_N \nabla_\xi \rho_0$ is of order $N$. This
explains the introduction of $\rho_0^\xi$.

The estimates just obtained prove the equicontinuity of the sequence
$\{\rho_t^{N,\xi} : N\ge 1\}$ on each compact set of $\bb R_+\times
\bb R$. Since every limit point is a weak solution of the heat
equation, by uniqueness of weak solutions, $\rho_t^\xi$ converges
uniformly on compact sets to the solution of \eqref{f8} with initial
condition $\gamma \rho_0$.

Since $\rho_0^\xi - \gamma \rho_0$ converges uniformly to 0, by the maximum
principle, $\rho_t^N$ converges uniformly on compact sets to the
solution of \eqref{f8}. This concludes the proof of the lemma.
\end{proof}

Let $h_t(x) = \xi_x (\nabla_N \rho_t^N)(x)$. A simple computation
shows that $(d/dt) h_t(x) = \xi_x (\Delta_N h_t)(x)$. Hence, $h_t$
satisfies a maximum principle and
\begin{equation}
\label{f28}
\sup_{t\ge 0} \, \sup_{x\in\bb Z} |(\nabla_N \rho_t)(x)| \;\le\;
\sup_{x\in\bb Z} |(\nabla_N \rho_0)(x)|
\end{equation}
because $\nabla_N \rho_t^N$ is absolutely bounded above and below by
$C(\epsilon) |h_t|$.

\medskip
\noindent{\bf Acknowledgments.} The authors would like to thank the
valuable discussions with A. Faggionato and S. R. S. Varadhan.

\end{document}